\documentclass[11pt]{article}
\usepackage{amssymb,latexsym,amsmath,caption,subcaption,float}
\usepackage{epstopdf}
\usepackage{graphicx,times}
\usepackage{mathtools}
\usepackage{varwidth}
\usepackage{breqn,caption,subcaption}
\usepackage[T1]{fontenc}
\newtheorem{The}{Theorem}[section]
\newtheorem{Def}{Definition}[section]

\newtheorem{Ex}{Example}[section]
\begin{document}
\begin{center}
{\Large \bf {A New family of methods for solving  delay differential equations}} \vskip 1cm
{\Large   {\bf Yogita Mahatekar, Pallavi S. Scindia}}\\

\textit{ Department of Mathematics, College of Engineering Pune, Pune - 411005, India, \\  yvs.maths@coep.ac.in,\, pus.maths@coep.ac.in}\\
\vskip1cm
{\large \textbf{Abstract}}\\
\end{center}
 In the present paper, we introduce a new family of $
 \theta-$methods for solving delay differential equations. New methods are developed using a combination of decomposition technique viz. new iterative method proposed by Daftardar Gejji and Jafari and  existing implicit numerical methods. Using Butcher tableau, we observed that new methods are non Runge-Kutta methods. Further, convergence  of new methods is investigated along with its stability analysis. Applications to variety of problems indicates that the proposed family of methods is more efficient than existing methods.
 \section{Introduction}
 A delay differential equation (DDE) is a differential equation in which state function is given in terms of value of the function at some previous times. Introduction of delay term in modelling allows better representation of real life phenomenon and enriches its dynamics. Due to presence of delay terms in the model, Delay differential equations (DDEs) are infinite dimentional and hence are difficult to analyse. Hence now a days, solving delay differential equations is an important area of research. Every DDE cannot be integrated analytically and hence there is a need to be dependant on numerical methods to solve DDEs. To develop efficient, stable and accurate numerical algorithms is primarily important task of research.
 \par DDEs are receiving increasing importance in many areas of science and engineering like biological processes, population growth and decay models, epidemiology, physiology, neural networks etc \cite{rihan2014delay}, \cite{bocharov2000numerical}, \cite{baker1998modelling}.  The classical numerical methods such as Eulers method, trapezoidal method, Runge-Kutta methods are discussed in \cite{butcher2000numerical}. Enright and Hu \cite{enright1995interpolating} have developed a tool to solve DDEs with vanishing delay with iteration and interpolation technique. Karoui and Vaillancourt \cite{karoui1995numerical} presneted a SYSDEL code to solve DDEs using Runge kutta methods of desired convergence order.
 In \cite{shampine2000solving}, a new MATLAB program dde23 has been developed to solve wide range of DDEs with constant delays. A new Adomian decomposition method is given in \cite{evans2005adomian} to solve DDEs. Further in 
\cite{ishak2008two}, two point predictor-corrector block method for solving DDEs is described. 
  Recently New iterative method (NIM) developed by Daftardar-Gejji and Jafari is used to develop many efficient numerical methods to solve fractional differential equations \cite{daftardar2015solving}, partial differential equations \cite{daftardar2008solving}, boundary value problems \cite{daftardar2010solving}. In \cite{aboodh2018solving}, a solution of DDE is acheived by using Aboodh transformation method. Recently in \cite{dehestani2020numerical} , a new numerical technique for solving fractional generalised pantograph-delay differential equations by using fractional order hybrid Bessel functions is developed.    
 \par 
 In the present article, we use the powerful technique of NIM
   to generate new efficient numerical tools to solve DDEs which are reducible to solve ODEs. 
\par The paper is organized as follows.
In section \ref{pre}, important preliminaries like Delay differential equations (DDEs), New iterative method (NIM) etc. are reviewed. In next section \ref{newmethods}, we developed a new family of numerical methods to solve DDEs and system of DDEs. In section \ref{err}, error analysis of newly proposed methods is done along with its stability analysis in section \ref{stab}. In section \ref{ex}, some illustrative examples are solved to check the accuracy of new methods practically. In last section \ref{conclusion}, some important observations are made on the basis of  theoretical stability analysis and error analysis of newly proposed methods.
\section{ Preliminaries}\label{pre}
\subsection{Delay differential equations}
Consider a general form of an initial value problem (IVP) representing a time delay differential equation:
\begin{equation} y' = f(t,y(t),y(t-\tau));\,t \in [0,T],\label{dde}\end{equation}
\begin{equation}
	y(t) = \phi(t),\,\, -\tau \leq t \leq 0, \label{em9}
\end{equation} where $\phi(t):[-\tau,0] \rightarrow \mathbb{R}^n;\, n \in \mathbb{N}$ is a real valued function which represents histroy of the solution in the past.
Existence and uniqueness conditions for the solution of DDEs are given in  \cite{suli2003introduction}. For existence and uniqueness of the solution of the IVP (\ref{dde}-\ref{em9}), we assume that $f$ is  non-linear, bounded and continuous real valued operator defined on
$[0,T]\times \mathbb{R}^n\times C^1(\mathbb{R},\mathbb{R}^n)$ to $\mathbb{R}^n$ and fulfils  Lipschitz conditions with
respect to the second and third arguments:
\begin{equation}|f(t,x_{1},u)-f(t,x_{2},u)|\leq L_{1} |x_{1}-x_{2}|\label{em10},\end{equation}\begin{equation}|f(t,x,u_{1})-f(t,x,u_{2})|\leq L_{2} |u_{1}-u_{2}|\label{em11},\end{equation} where $L_{1}$ and $L_{2}$ are positive constants.
\subsection{Approximation to the delay term $y(t-\tau)$}\label{approxdelay}
The approximation to the delay term  $y(t_n-\tau)$ is denoted by $\nu_{n}$ \cite{daftardar2015solving}.\\
When $\tau$ is constant $(t_{n}-\tau)$ may not be a grid point
$t_{n}$ for any n. Suppose $(m+\delta)h=\tau, m\in\mathbb{N}$ and $0\leq\delta<1.$\\
If $\delta=0,\,mh=\tau$ and $y(t_{n}-\tau)$ is approximated as
$$
y(t_{n}-\tau)\approx \nu_{n}=\left\{\begin{array}{ccc}
	y_{n-m}&\mbox{ if}&n\geq m;\\
	\phi({t_{n}-\tau})&\mbox{if}&n< m.\\
\end{array}\right.
$$
\subsection{New iterative method (NIM)}\label{eulernim}
\par Daftardar-Gejji and Jafari \cite{daftardar2006iterative} have proposed a
new iterative method (NIM) for solving linear/non-linear functional
equations of the form
\begin{equation}
u=f + L(u)+ N(u)\label{em7}
\end{equation}
where $f$ is a known function, $L$ is a linear operator and $N$ a non linear operator. This decomposition technique is used by many researchers for solving variety of problems such as fractional differential equations \cite{daftardar2015solving}, boundary value problems \cite{daftardar2010solving}, system of non-linear functional equations \cite{bhalekar2012solving} etc. Recently NIM is used to develop efficient numerical algorithms to solve ordinary differential equations and time delay fractional differential equations too \cite{sukale2017new},\cite{jhinga2019new}.
\par In this method, We assume that eq.(\ref{em7}) has a series solution of the form $u=\displaystyle\sum^{\infty}_{i=0}u_{i}.$ Since $L$ is a linear operator, we have $L\left(\displaystyle\sum^{\infty}_{i=0}u_{i}\right)=\displaystyle\sum^{\infty}_{i=0}L(u_{i})$ and non linear operator $N$ is decomposed as
$$N(u)=N(u_{0})+[N(u_{0}+u_{1})-N(u_{0})]+[N(u_{0}+u_{1}+u_{2})-N(u_{0}+u_{1})]+\dots$$
Let $G_{0}=N(u_{0})\,\text{and}\,
G_{i}=N\left(\displaystyle\sum_{n=0}^{i}u_{n}\right)-N\left(\displaystyle\sum_{n=0}^{i-1}u_{n}\right), i=1,2,3,\dots$. \\ Observe that
$N(u) =\displaystyle\sum_{i=0}^{\infty}G_{i}$.\\
Putting series solution $u=\displaystyle\sum^{\infty}_{i=0}u_{i}$ in eq.(\ref{em7}), we get $$\displaystyle\sum^{\infty}_{i=0}u_{i}=f+\displaystyle\sum^{\infty}_{i=0}L(u_{i})+\displaystyle\sum_{i=0}^{\infty}G_{i}.$$
Taking $u_{0}=f$, and $u_{n}=L(u_{n-1})+G_{n-1}, n=1,2,3,\dots$
\\Note that $$u=u_{0}+u_{1}+u_{2}+\dots =
f+L(u_0)+N(u_{0})+L(u_1)+[N(u_{0}+u_{1})-N(u_{0})]+\dots=f+L(u)+N(u).$$ Hence $u$
satisfies the functional eq.(\ref{em7}). $k$-term NIM solution is given by $u=\displaystyle\sum^{k-1}_{i=0}u_{i}.$
\section{New methods to solve DDEs}\label{newmethods}
In this section, we represent a family of new numerical methods based on NIM and implicit numerical methods to solve delay differential equations.
In \cite{sukale2017new}, we developed a new algorithm to solve DDEs and ODEs by improving existing trapezoidal rule of integration. In present paper, we generalize this task by improving general $\theta-$methods to solve delay differential equations as follows. 
\par For solving eq.(\ref{dde}-\ref{em9}) on $[0,T]$, consider the uniform grid $t_n=nh,\,n=-m,-m+1,....,-1,0,1,....N,$ where $m$ and $N$ are integers such that $N=T/h$ and $m=\tau/h.$
We integrate eq.(\ref{dde}) from the node $t_{n}$ to $t_{n+1}$ on both sides, which gives us \begin{equation}y(t_{n+1})=y(t_{n})+\displaystyle\int^{t_{n+1}}_{t_{n}}f(t,y(t),y(t-\tau))\,\,dt.\label{em12}\end{equation}
Now we approximate integration on R.H.S. in above equation as;
\begin{equation}\displaystyle\int^{t_{n+1}}_{t_n}f(t,y(t),y(t-\tau))dt\approx h \left[(1-\theta)f(t_n,y(t_n),\nu_n)+\theta f(t_{n+1},y(t_{n+1}),\nu_{n+1})\right],\label{em4} \end{equation}
where $\theta$ is a parameter which lies in the closed interval $[0,1]$ and $\nu_n$ is approximation to delay term at the node $t=t_n.$  Applying the approximation from eq.(\ref{em4}) in eq.(\ref{em12}), we get 
\begin{equation}y_{n+1}=y_n+h\left[(1-\theta)f(t_n,y_n,\nu_n)+\theta f(t_{n+1},y_{n+1},\nu_{n+1})\right],\,n=0,1,\hdots , N-1.\label{em5} \end{equation}
This is called as a family of $\theta-$ methods to solve DDEs.\\
In particular, when $\theta=1,$ eq.(\ref{em5}) reduces to 
\begin{equation}y_{n+1}=y_n+hf(t_{n+1},y_{n+1},\nu_{n+1}),\,n=0,1,\hdots , N-1.\label{em6} \end{equation}
This  eq.(\ref{em6}) is referred as Implicit Euler's method.
Note that when $\theta=0,$ eq.(\ref{em5}) yields Explicit Euler's method and when $\theta =\displaystyle\frac{1}{2}$ we get implicit Trapezoidal rule of integration to solve DDEs. 
\subsection{New Reults}
To generate a new family of $\theta-$methids for solving DDEs, we note that eq.(\ref{em5}) is of the form $u=f+N(u),$ and hence NIM can be employed as follows:
\par We write, 
	$$u=y_{n+1},$$
	 $$f=y_n+h(1-\theta)f(t_n,y_n,\nu_n),$$ 
	$$ N(u)=h\theta f(t_{n+1},y_{n+1},\nu_{n+1}).$$

For simplicity we denote, $f_n=f(t_n,y_n,\nu_n).$ Three term NIM solution of eq.(\ref{em5}) gives us 
\begin{eqnarray*}
u&=&u_0+u_1+u_2\\
&=&u_0+N(u_0)+N(u_0+u_1)-N(u_0)\\
&=&u_0+N(u_0+u_1)\\
&=&u_0+N(u_0+N(u_0)).\end{eqnarray*}
That is \\
$y_{n+1}=y_n+h(1-\theta)f_n+N(u_0+u_1)$\\
$=y_n+h(1-\theta)f_n+h\theta f(t_{n+1},u_0+u_1,\nu_{n+1})$\\
Therefore,\begin{equation}y_{n+1}=y_n+h(1-\theta)f_n+h\theta f(t_{n+1},y_n+h(1-\theta)f_n+h\theta f(t_{n+1},y_n+h(1-\theta)f_n,\nu_{n+1}),\nu_{n+1})\label{newfamily}\end{equation}
Eq.(\ref{newfamily}) represents new family of $\theta-$methods to solve DDEs which can be expressed in the following more simple form as follows:\\
\begin{equation}
\begin{rcases}
	\begin{array}{r@{}l}
		k_1=f(t_n,y_n,\nu_n)\\
	k_2=f(t_{n+1},y_n+h(1-\theta)k_1,\nu_{n+1})
	\end{array}\\
		\begin{array}{r@{}l}
		k_3=f(t_{n+1},y_n+h(1-\theta)k_1+h\theta k_2,\nu_{n+1})\\
	\text{where,\,\,}	y_{n+1}=y_n+h(1-\theta)k_1+h\theta k_3.
	\end{array}\\
\end{rcases}\begin{array}{l}
\end{array}\label{theta}\end{equation}
{\bf Case 1:} When $\theta=1,$\\
\begin{equation}
	\begin{rcases}
		\begin{array}{r@{}l}
		k_1=f(t_n,y_n,\nu_n)\\
			k_2=f(t_{n+1},y_n,\nu_{n+1})
		\end{array}\\
		\begin{array}{r@{}l}
			k_3=f(t_{n+1},y_n+hk_2,\nu_{n+1})\\
			\text{where,\,\,}	y_{n+1}=y_n+h k_3.
		\end{array}\\
	\end{rcases}\begin{array}{l}
	\end{array}\label{theta1}
\end{equation}Eqs.(\ref{theta1}) represents new improved implicit Eulers method to solve DDEs.\\
{\bf Case 2:} When $\theta=\frac{1}{2},$\\
\begin{equation}
	\begin{rcases}
		\begin{array}{r@{}l}
			k_1=f_n=f(t_n,y_n,\nu_n)\\
		k_2=f(t_{n+1},y_n+\frac{hk_1}{2},\nu_{n+1})
		\end{array}\\
		\begin{array}{r@{}l}
			k_3=f(t_{n+1},y_n+\frac{hk_1}{2}+\frac{hk_2}{2},\nu_{n+1})\\
			\text{where,\,\,}	y_{n+1}=y_n+\frac{hk_1}{2}+\frac{hk_3}{2}.,
		\end{array}\\
	\end{rcases}\begin{array}{l}
	\end{array}\label{theta12}
\end{equation}
Eqs.(\ref{theta12}) represents improved trapezoidal rule to solve DDEs which is developed in \cite{sukale2017new}.\\
{\bf Case 3:} When $\theta=0,$\\
\begin{equation}
	\begin{rcases}
		\begin{array}{r@{}l}
			k_1=f_n=f(t_n,y_n,\nu_n)\\
			k_2=f(t_{n+1},y_n+hk_1,\nu_{n+1})
		\end{array}\\
		\begin{array}{r@{}l}
			k_3=f(t_{n+1},y_n+hk_1,\nu_{n+1})\\
			\text{where,\,\,}	y_{n+1}=y_n+hk_1.,
		\end{array}\\
	\end{rcases}\begin{array}{l}
	\end{array}\label{theta0}
\end{equation}

Eqs.(\ref{theta0}) represents a method which is usual explicit Eulers method to solve DDEs.\\
{\bf Case 4:} When $\theta=\frac{3}{4},$\\
\begin{equation}
	\begin{rcases}
		\begin{array}{r@{}l}
			k_1=f_n=f(t_n,y_n,\nu_n)\\
			k_2=f(t_{n+1},y_n+\frac{hk_1}{4},\nu_{n+1})
		\end{array}\\
		\begin{array}{r@{}l}
			k_3=f(t_{n+1},y_n+\frac{hk_1}{4}+\frac{3hk_2}{4},\nu_{n+1})\\
			\text{where,\,\,}y_{n+1}=y_n+\frac{hk_1}{4}+\frac{3hk_3}{4}.,
		\end{array}\\
	\end{rcases}\begin{array}{l}
	\end{array}\label{theta34}
\end{equation}

Eqs.(\ref{theta34}) is a new method to solve DDEs.
\par In above family of methods to solve DDEs, when delay term $y(t_n-\tau)$ and its approximation $\nu_n$ are zero then DDE (\ref{dde}) reduces to ODE without delay. In accordance with this, above family of methods gets reduced to the family of numerical methods for solving ODEs (without delay) which are developed in \cite{ababneh2019new}.
\subsection{Non-Runge Kutta methods}
General form of Runge-Kutta method is given by

$y_{n+1}=y_{n}+h\displaystyle\sum^{3}_{i=1}b_ik_i$ and 
\begin{equation}
	\begin{rcases}
		\begin{array}{r@{}l}
			k_1=f_n=f(t_n,y_n,\nu_n)\\
			k_2=f(t_{n}+c_2h,y_n+ha_{21}k_1,\nu_{n})
		\end{array}\\
		\begin{array}{r@{}l}
			k_3=f(t_{n}+c_3h,y_n+h (a_{31}k_1+a_{32}k_2),\nu_{n})\\
		\end{array}\\
	\end{rcases}\begin{array}{l}
	\end{array}\label{rungekutta}
\end{equation}
In newly proposed $\theta-$methods represented by eqs.(\ref{theta}), we have $b_{1}(\theta)=1-\theta,\,b_{2}(\theta)=0,\, b_{3}(\theta)=\theta,\, a_{21}=1-\theta,\, a_{31}=1-\theta,\, a_{32}=\theta,\, c_{2}=1,\, c_{3}=1.$ Therefore new Family of $\theta-$methods for solving DDEs can be stated in the form of Butcher tableu as follows.

\begin{center}
	\begin{tabular}{c|cccc}
		0& & & \\
		1&1-$\theta$& & \\
		1&1-$\theta$&$\theta$&\\  
		\hline\\
		&1-$\theta$&0&$\theta$&
	\end{tabular}
\end{center}
For a Runge Kutta method it is necessary to satisfy that $\displaystyle\sum^{i-1}_{j=1} a_{ij}=c_{i}\, \forall\, i=2,3\,\, (cf. $\cite{bellen2013numerical}). From the above tableau, for $i=2,\,$ $\displaystyle\sum^{1}_{j=1} a_{2j}=c_{2}$ if and only if $\theta=0.$ This shows that, newly proposed family of numerical methods for solving DDEs is different from Runge Kutta methods except for $\theta=0.$

 \subsection{New $\theta-$methods for solving a system of delay differential equations}\label{syst3}
 The numerical algorithm presented in above section  can be generalized for solving the following system of DDEs: \\
 $$y'_{1}(t)=f_{1}(t,\overline{y}(t),\overline{y}(t-\tau)),$$
 $$y'_{2}(t)=f_{2}(t,\overline{y}(t),\overline{y}(t-\tau)),$$
 $$\vdots$$
 $$y'_{m}(t)=f_{m}(t,\overline{y}(t),\overline{y}(t-\tau)),$$ with initial condition \begin{equation}\overline{y}(t)=(y_{1}(t),y_{2}(t),\hdots,y_{m}(t))=(\phi_{1}(t),\phi_{2}(t),...,\phi_{m}(t));-\tau\leq t\leq 0.\label{syst1}\end{equation}
 We let $\overline{y_n}$ be a vector of independent variables representing values of $(y_1,y_2,\hdots,y_m)$ at node $t=t_n$ and $\overline{\nu_n}$ is a vector approximation to $(y_{1}(t-\tau),y_{2}(t-\tau),\hdots,y_{m}(t-\tau))$ at $t=t_n.$
 We obtain a new family of $\theta-$methods to solve a system of DDEs as follows:\\
$k_{1,y_i}=f_{i}(t_n,\overline{y_n},\overline{\nu_n})$\\
$k_{2,y_i}=f_{i}(t_{n+1},\overline{y_n}+h(1-\theta)k_{1,y_i},\overline{\nu_{n+1}})$\\
$k_{3,y_i}=f_{i}(t_{n+1},\overline{y_n}+h(1-\theta)k_{1,y_i}+h\theta k_{2,y_i},\overline{\nu_{n+1}})$\\
Where,$$\overline{y_{n+1}}=\overline{y_n}+h(1-\theta)k_{1,y_i}+h\theta k_{3,y_i},\,i=1,2,\hdots,m.$$
  \section{Error analysis}\label{err}
  \begin{The}
  	The new family of $\theta-$methods given by eqs.(\ref{theta}) forms a second order numerical method  for $\theta=\frac{1}{2}$ and has a first order  convergence for any other value of $\theta \in [0,1].$\end{The}
  Proof: Using Taylor's series expansion,
  $$y(t_{n+1})=y(t_n+h)=y_n+hf_n+\frac{h^2}{2}(f_t+ff_y+f_\nu)+O(h^3).$$
  Now consider $\theta-$method and Taylors expansions applied in $k_2,\,k_3.$
  		\begin{eqnarray}
  			k_1&=&f(t_n,y_n,\nu_n)\label{k1}
  			\end{eqnarray}
  			\begin{eqnarray*}
  			k_2&=&f(t_{n+1},y_n+h(1-\theta)k_1,\nu_{n+1})\\
  			&=&f(t_{n}+h,y_n+h(1-\theta)k_1,\nu_{n}+h)\\
  			  		\end{eqnarray*}  		
  \begin{dmath}k_2=f_n+\left(hf_t+h(1-\theta)k_1f_y+hf_{\nu}\right)+\\
  \frac{1}{2}\left(hf_{tt}+h^2(1-\theta)^2k^2_1f_{yy}+h^2f_{\nu\nu}+2h^2(1-\theta)k_1f_{ty}+2h^2f_{t\nu}+2h^2(1-\theta)k_1f_{y\nu}\right)+O(h^3).\label{k2}\end{dmath}
  \begin{eqnarray*}
  	  	k_3&=&f(t_{n+1},y_n+h(1-\theta)k_1+h\theta k_2,\nu_{n+1})\\
  	&=&f(t_{n}+h,y_n+h(1-\theta)k_1+h\theta k_2,\nu_{n}+h)\\
  	\end{eqnarray*}
  \begin{dmath}
  	k_3=f_n+\left( hf_t+h(1-\theta)k_1f_y+h\theta k_2f_y+hf_{\nu}\right)+
  	\frac{1}{2}\left(h^2f_{tt}+(h^2(1-\theta)^2k^2_1+h^2\theta^2k^2_2+2h^2(1-\theta)\theta k_1k_2)f_{yy}\\+h^2f_{\nu\nu}+(2h^2(1-\theta)k_1+2h^2\theta k_2)f_{yt}+2h^2f_{t\nu}+2(h^2(1-\theta)k_1+h^2\theta k_2)f_{y \nu}\right)+O(h^3).\label{k3}
  \end{dmath}
  Putting (\ref{k1}), (\ref{k2}), (\ref{k3}) in $\theta-$method given by eqs.(\ref{theta}), error $e_{n+1}$ of the method is given by
  \begin{dmath}
  	e_{n+1}=y(t_{n+1})-y_{n+1}\\
  	=y_n+hf_n+\frac{h^2}{2}f_t+\frac{h^2ff_y}{2}+\frac{h^2f_{\nu}}{2}+O(h^3)-\left(y_n+h(1-\theta)f_n+h\theta f_n\\
  	h^2\theta f_t+h^2 \theta (1-\theta)k_1 f_y+h^2 \theta^2 k_2 f_y+h^2 \theta f_{\nu}\right)+O(h^3).
  \end{dmath}
Clearly, for $\theta=0,1,3/4$ method has a linear convergence and for $\theta=1/2$ its order of convegence is quadratic. 
\subsection{Error analysis of new improved implicit Euler's method for solving DDEs
 } Refer to Eqs.(\ref{theta1}), For $\theta=1,$ we get new improved implicit Eulers method to solve DDEs, which is given by the formula: 
\begin{equation}
y_{n+1}=y_n+hf[t_{n+1},y_n+hf(t_{n+1},y_n,\nu_{n+1}),\nu_{n+1}].\label{em22}\end{equation}
By inserting analytical solution in the above equation, we get the truncation errorv $T_n$ as 
\begin{equation}
\displaystyle\frac{y(t_{n+1})-y(t_n)}{h}-f[t_{n+1},y(t_n)+hf(t_{n+1},y(t_n),\nu_{n+1}),\nu_{n+1}]=T_n.\label{em23}\end{equation} 
and by eq.(\ref{em23}), we get 
\begin{equation}
\displaystyle\frac{y_{n+1}-y_n}{h}-f[t_{n+1},y_n+hf(t_{n+1},y_n,\nu_{n+1}),\nu_{n+1}]=0.\label{em24}\end{equation} 
Subtracting eq.$(\ref{em23})$ and eq.$(\ref{em24})$, we get 
\begin{equation}
T_n=\displaystyle\frac{e_{n+1}-e_n}{h}-f[t_{n+1},y(t_n)+hf(t_{n+1},y(t_n),\nu_{n+1}),\nu_{n+1}]+f[t_{n+1},y_n+hf(t_{n+1},y_n,\nu_{n+1}),\nu_{n+1}].\label{em25}\end{equation}
This implies that,
\begin{dmath*}
hT_n=e_{n+1}-e_n-h\left(f[t_{n+1},y(t_n)+hf(t_{n+1},y(t_n),\nu_{n+1}),\nu_{n+1}]\\-f[t_{n+1},y_n+hf(t_{n+1},y_n,\nu_{n+1}),\nu_{n+1}]\right).\label{em26}\end{dmath*}
Therefore,
\begin{dmath*}
	|e_{n+1}|\leq |e_n|+hL_1|(y(t_n)+hf(t_{n+1},y(t_n),\nu_{n+1})-(y_n+hf(t_{n+1},y_n,\nu_{n+1})|+hT_n\\	
	\leq |e_n|+hL_1 |e_n|+h^2L_1^2 |e_n|+|T_n|h\\
	\leq (1+hL_1+h^2L_1^2)|e_n|+|T_n|h;\,\,n=0,1,\hdots N.
\end{dmath*}
Let $T=\displaystyle\max_{0\leq n \leq (N-1)}|T_n|$
Therefore, 
$$|e_{n+1}|\leq (1+hL+h^2L^2)|e_n|+Th$$
Now by induction,
\begin{dmath*}
	|e_n| \leq (1+hL+h^2L^2)^n|e_0|+\left(\frac{(1+hL+h^2L^2)^{n-1}-1}{(1+hL+h^2L^2)-1}\right)Th\\ \leq (1+hL+h^2L^2)^n|e_0|+\left(\frac{e^{nhL}-1}{(Lh+1)Lh}\right)Th\\
	\leq e^{nhL}|e_0|+\left(\frac{e^{nhL}-1}{L}\right)T\\
	\leq e^{(t_n-t_0)L}|e_0|+\left(\frac{e^{(t_n-t_0)L}-1}{L}\right)T
\end{dmath*}
Since $nh=t_n-t_0.$
Noting that, Truncation error $T \leq \frac{1}{2}hy''(\zeta)$ and if $y'' \leq M_2$ then $T \leq \frac{hM_2}{2}.$ Therefore,

\begin{dmath*}|e_n| \leq e^{(t_n-t_0)L_1}|e_0|+\left(\displaystyle\frac{(e^{(t_n-t_0)L_1}-1)hM_2}{2L_1}\right)\\ \leq \left(\displaystyle\frac{(e^{(t_n-t_0)L_1}-1)hM_2}{2L_1}\right), \text{noting that $|e_0|$ is zero.}
\end{dmath*}hence as $h \rightarrow 0,$ $|e_n|\rightarrow 0.$ Here $L_1$ is Lipschitz constant as given in eq.(\ref{em10}). Hence method is convergent.
Similary, for $\theta=\frac{1}{2},$ we proved the convergence of the new improved trapezoidal rule for solving DDEs in \cite{sukale2017new} and for other values of $\theta$ convergenrce can be proved on similar lines.

\section{Stability Analysis}\label{stab}
\begin{Def}
The numerical method for solving IVP eq.(\ref{dde}-\ref{em9}) is said to be zero-stable if small perturbation in the initial condition of IVP do not cause the numerical approximation to diverge from the exact solution, provided the exact solution of the IVP is bounded. 
\end{Def}
\par Consider the IVP eq.(\ref{dde}-\ref{em9}) and let $\epsilon{y}(0)=\epsilon{y_{0}}$ be the new initial value (perturbed initial condition) obtained by making a small
change in $y(0)=y_{0}.$
\begin{The}
Let $y_n$ be the solution obtained by new improved implicit Euler's  method (Case 1 in section\ref{newmethods}) for solving DDE at the node $t_n$ with the initial condition $y(0)=y_0$ and let $\epsilon y_n$ be the solution obtained by the same numerical method with perturbed initial condition $\epsilon{y_{0}}=y_0+\epsilon_0;\,\epsilon_0 >0.$ We assume that $f(t,y)$ satisfies Lipschitz condition with respect to second variable and third variable with Lipschitz constant $L_1,\,L_2$ then $\exists$ positive constants $k$ and $\epsilon_1$ such that $|y_n-\epsilon y_n| \leq k \epsilon,\,\, \forall nh \leq T,\,\, h \in (0, \epsilon_1)\, $ whenever $|\epsilon_0|\leq \epsilon.$
\end{The}
Proof: We prove the result by induction.
 \begin{dmath*}|y_n-\epsilon y_n| = |y_{n-1}+hf\left(t_n,y_{n-1}+hf\left(t_n,y_{n-1},\nu_n\right),\nu_n\right)-\epsilon y_{n-1}-hf\left(t_n,\epsilon y_{n-1}+hf\left(t_n,\epsilon y_{n-1},\nu_n\right),\nu_n\right)|\\
 	\leq |y_{n-1}-\epsilon y_{n-1}|+hL_1|y_{n-1}+hf\left(t_n,y_{n-1},\nu_n\right)-\epsilon y_{n-1}-hf\left(t_n,\epsilon y_{n-1},\nu_n\right)|\\
 		\leq |y_{n-1}-\epsilon y_{n-1}|+hL_1|y_{n-1}-\epsilon y_{n-1}|+h^2L_1^2|y_{n-1}-\epsilon y_{n-1}|\\
 		=(1+hL_1+h^2L_1^2)|y_{n-1}-\epsilon y_{n-1}|
\end{dmath*}  
Therefore by induction,
\begin{dmath*}|y_n-\epsilon y_n| \leq (1+hL_1+h^2L_1^2)^n|y_{0}-\epsilon y_0|\\
	\leq e^{nhL_1} |\epsilon_0|\\
	=e^{T_{M}L_1} \epsilon = k \epsilon,\,\, \text{where}\, k=e^{T_{M}L_1}>0.
\end{dmath*}
This proves that the new method for ($\theta=1$) is stable.
\section{Illustrative examples}\label{ex}
To demonstrate applicability of newly proposed methods, We present here some illustrative examples which are solved using Mathematica 12.
\begin{Ex}
	Consider the  delay logistic differential equation 	\begin{eqnarray}
		y'(t) &=& 0.3 y(t)(1-y(t-1)),\, y(t \leq 0) = 0.1. \label{ex.2}
	\end{eqnarray}
\end{Ex}
\begin{figure}[H]
	\begin{subfigure}[b]{.5\linewidth}
		\centering
		\includegraphics[scale=.5]{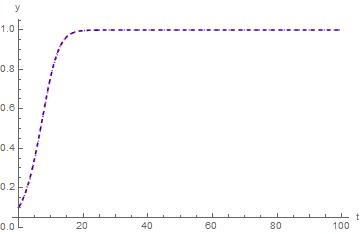}
		\caption{Solution of (\ref{ex.2}) by new  method ($\theta=1$)}
	\end{subfigure}%
	\begin{subfigure}[b]{.5\linewidth}
		\centering
		\includegraphics[scale=.5]{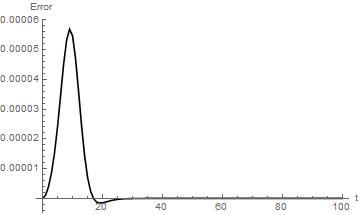}
		\caption{Error in solution of (\ref{ex.2}) when solved using  new  method ($\theta=1$)}
	\end{subfigure} 	  
	\caption{Dashed graph: Solution by new  method ($\theta=1$), Dotted graph:  Exact solution. It is noted that in Fig.\ref{fig:3} (a) exact solution and approximate solution by the  new  method ($\theta=1$) overlaps on each other.} Step length is taken as $h=0.01.$\label{fig:3}
\end{figure}
\begin{Ex}
	Consider the  differential equation without delay \cite{bui2010explicit}
	\begin{eqnarray}
	y'(t) &=& 2-e^{-4t}-2y,\, y(0) = 1,\, 0 \leq t \leq 10. \label{ex.1}
	\end{eqnarray}
\end{Ex}
Exact solution of the differential equation (\ref{ex.1}) is $y(t)=1+\displaystyle\frac{e^{-4t}-e^{-2t}}{2}.$
 \begin{figure}[H]
 	\begin{subfigure}[b]{.5\linewidth}
 		\centering
 		\includegraphics[scale=.5]{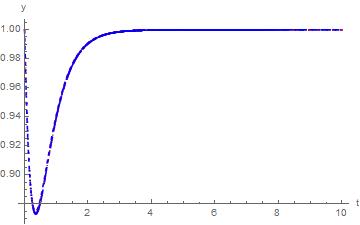}
 		\caption{Solution of (\ref{ex.1}) by new method ($\theta=1$)}
 	\end{subfigure}%
 	\begin{subfigure}[b]{.5\linewidth}
 		\centering
 		\includegraphics[scale=.5]{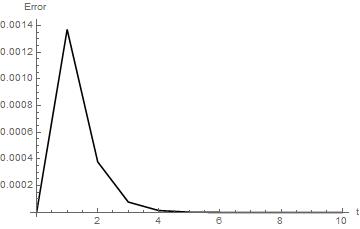}
 		\caption{Error in solution of (\ref{ex.1}) when solved using $\theta$ method ($\theta=1$)}
 	\end{subfigure} 	  
 	\caption{Dashed graph: Solution by new method, Dotted graph:  Exact solution. It is noted that in Fig.\ref{fig:2} (a) exact solution and approximate solution by the new method overlaps on each other.} Step length is taken as $h=0.01.$\label{fig:2}
 \end{figure}
The  error in the new method is shown in Fig.\ref{fig:2} (b).
Thus, the new method is accurate. In Table \ref{t1}, error in implicit backward Euler's method $E_{1}$  and error in new method ($\theta=1$) $E_{2}$ (obtained by taking 3-term NIM solution) while solving \eqref{ex.1} are compared for various values of $t$. It is noteworthy that $E_{2}$ is always smaller than $E_{1}$ in all cases. So the  new  method is more accurate than implicit backward Euler's method.\\

  \begin{table}[H]
  	\begin{tabular}[H]{|c|c|c|c|c|c|c|c|}
  		\hline
  		$t_n$ & $S_1$  & $S_2$ & $S_3$ & S & $e_{1}$ & $e_{2}$ & $e_{3}$ \\
  		\hline
  		$0$&1
  		&1 & $1$ &$ 1$ & $ 0$&$ 0$&$ 0$ \\
  		\cline{1-8}	$0.01$& $0.99016$
  		&$0.9902$ & $0.990196$ &$ 0.990295$ & $0.000099$&$0.0000954$&$ 0.0000994$ \\
  		\hline
  	$0.02$& $0.980969$
  	&$0.980976$ & $0.980969$ &$ 0.981163$ & $0.000194$&$0.000187$&$ 0.0001948$ \\
  		\cline{1-8}		$0.03$& $0.972292$
  		&$0.97230$ & $0.972292$ &$ 0.972578$ & $0.000286$&$0.000275$&$ 0.0002863$ \\
  		\hline
  			$0.04$& $0.96414$
  			&$0.9641$ & $0.96414$ &$ 0.964514$ & $0.0003739$&$0.00036$&$ 0.0003741$ \\
  		\cline{1-8}		$0.05$& $0.95648$
  		&$0.9565$ & $0.956488$ &$ 0.956947$ & $0.000458$&$0.0004413$&$ 0.0004584$ \\
  		\hline
  				$0.06$& $0.949315$
  				&$0.949334$ & $0.949315$ &$ 0.949854$ & $0.000538$&$0.0005194$&$ 0.0005391$ \\
  		\cline{1-8}		$0.1$& $0.924966$
  		&$0.924993$ & $0.924966$ &$ 0.925795$ & $0.0008283$&$0.00080$&$ 0.0008289$ \\
  		\hline
  				$0.2$& $0.980969$
  				&$0.980976$ & $0.980969$ &$ 0.981163$ & $0.000194$&$0.000187$&$ 0.0001948$ \\
  		\cline{1-8}		$0.3$& $0.8745$
  		&$0.874563$ & $0.874528$ &$ 0.87619$ & $0.0016622$&$0.001628$&$ 0.001663$ \\
  		\hline
  				$0.4$& $0.87447$
  				&$0.874498$ & $ 0.87447$ &$0.87628$ & $0.001813$&$0.001785$&$ 0.001814$ \\
  		\cline{1-8}		$0.5$& $0.881873$
  		&$0.881893$ & $0.881872$ &$ 0.883728$ & $0.001855$&$0.001834$&$ 0.0018554$ \\
  		\hline
  	\end{tabular}
  	\caption{ Ex.(8.1)\\
  		$S_1:$ Solution by backward Euler's method\\
  		$S_2:$ Solution by new method ($\theta=1$) obtained by taking 3-term NIM solution\\
  		$S_3:$ Solution by new method ($\theta=1$) obtained by taking 4-term NIM solution\\
  		$S:$ Exact solution\\
  			$e_1:$ Error in solution by backward Euler's method\\
  			$e_2:$ Error in solution by new method ($\theta=1$) obtained by taking 3-term NIM solution\\
  			$e_3:$ Error in solution by new method ($\theta=1$) obtained by taking 4-term NIM solution\\ }
  	\label{t1}
  \end{table}
{\bf Observation:} In this example, it is observed that (4-term NIM solution) new method ($\theta=1$) method and backward Euler's method gives same error and error in these two methods is greater than (3-term NIM solution) new method ($\theta=1$). Hence, new method with three term NIM solution gives better accuracy than implicit backward euler method and new method with 4-term NIM solution.
\subsection{R$\ddot{o}$ssler System with delay}
\begin{figure}[H]
	\begin{subfigure}[b]{.5\linewidth}
		\centering
		\includegraphics[scale=.5]{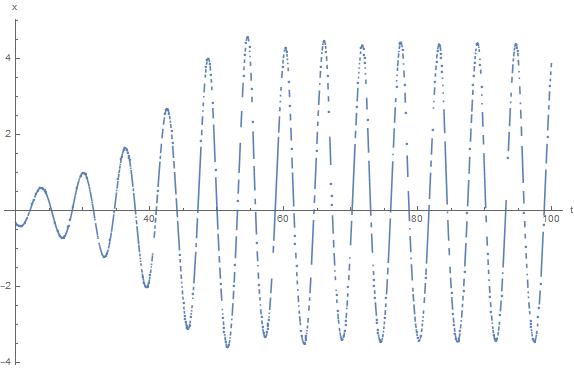}
		\caption{ x-waveform of the system \eqref{syst2}  }
	\end{subfigure}
	\begin{subfigure}[b]{.5\linewidth}
		\centering
		\includegraphics[scale=.5]{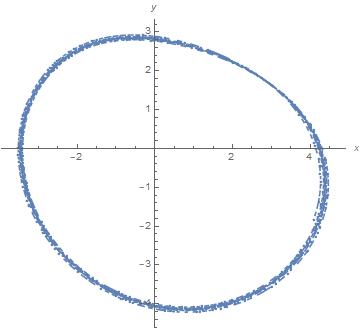}
		\caption{xy-phase diagram of the system  \eqref{syst2}}
	\end{subfigure} 	
	\caption{c=2.3 in system\eqref{syst2}}\label{fig:ros1}
\end{figure}

\begin{figure}[H]
	\begin{subfigure}[b]{.5\linewidth}
		\centering
		\includegraphics[scale=.5]{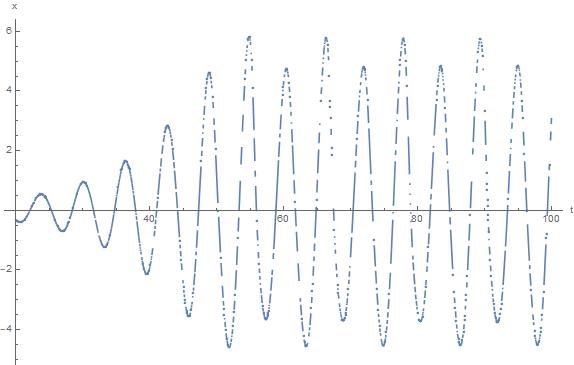}
		\caption{ x-waveform of the system \eqref{syst2}  }
	\end{subfigure}
	\begin{subfigure}[b]{.5\linewidth}
		\centering
		\includegraphics[scale=.5]{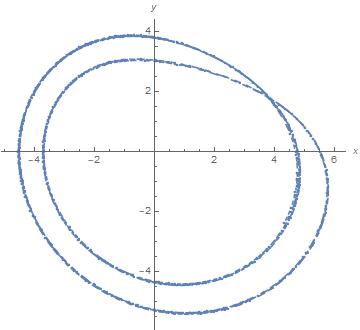}
		\caption{xy-phase diagram of the system  \eqref{syst2}}
	\end{subfigure} 	
	\caption{c=2.9 in system\eqref{syst2}}\label{fig:ros2}
\end{figure}

\begin{figure}[H]
	\begin{subfigure}[b]{.5\linewidth}
		\centering
		\includegraphics[scale=.5]{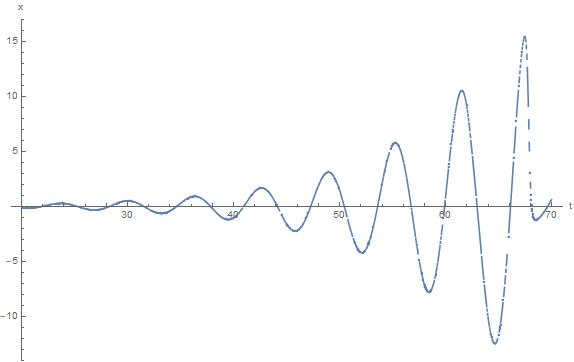}
		\caption{ x-waveform of the system \eqref{syst2}  }
	\end{subfigure}
	\begin{subfigure}[b]{.5\linewidth}
		\centering
		\includegraphics[scale=.5]{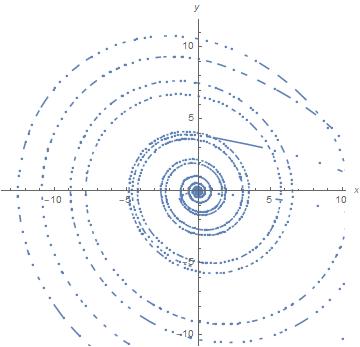}
		\caption{xy-phase diagram of the system  \eqref{syst2}}
	\end{subfigure} 	
	\caption{c=7.9 in system\eqref{syst2}}\label{fig:ros3}
\end{figure}

Consider the R$\ddot{o}ssler$ system  \cite{ibrahim2018chaotic} with delay given by the following system of differential equations.
\begin{eqnarray}
\dot{x}&&= -y(t)-z(t),\nonumber \\\dot{y}&&=x(t)+ay(t-1),\nonumber\\ \dot{z}&&=b+z(t)(x(t)-c).\label{syst2}
\end{eqnarray}
Let $a=b=0.2, \,\, x(0)=y(0)=z(0)=0.0001.$ Note that $c$ is a control parameter.  Solving above system by new method for ODE we obtain the x-waveforms  and x-y phase potraits which are depicted in Figs.\ref{fig:ros1}(a)-(b), Figs.\ref{fig:ros2}(a)-(b) and in  Figs.\ref{fig:ros3}(a)-(b) for $c=2.3,\,c=2.9,\,c=7.9$ respectively. \\

\section{Conclusions}\label{conclusion}
In the present work, a new family of methods for solving delay differential equations (DDEs) has been proposed which are reducible to solve ordinary differential equations too (without delay). Newly proposed methods are then compared with existing methods with respect to stability and accuracy. New methods formed found to be stable and accurate. Further error analysis and stability analysis of new methods is carried out. Numerous illustrative examples are solved using Mathematica 12 to demonstrate the efficiency of the method. It is observed that, new methods are non-Runge Kutta methods and are more accurate than existing numerical methods for solving DDEs. 
\vskip 1cm
\bibliography{reference1_xeriv}

\begin{thebibliography}{10}

\bibitem{ababneh2019new}
O.~Y. Ababneh.
\newblock New numerical methods for solving differential equations, 2019.

\bibitem{aboodh2018solving}
K.~Aboodh, R.~Farah, I.~Almardy, and A.~Osman.
\newblock Solving delay differential equations by aboodh transformation method.
\newblock {\em International Journal of Applied Mathematics \& Statistical
  Sciences}, 7(2):55--64, 2018.

\bibitem{baker1998modelling}
C.~Baker, G.~Bocharov, C.~Paul, and F.~Rihan.
\newblock Modelling and analysis of time-lags in some basic patterns of cell
  proliferation.
\newblock {\em Journal of mathematical biology}, 37(4):341--371, 1998.

\bibitem{bellen2013numerical}
A.~Bellen and M.~Zennaro.
\newblock {\em Numerical methods for delay differential equations}.
\newblock Oxford University Press, 2013.

\bibitem{bhalekar2012solving}
S.~Bhalekar and V.~Daftardar-Gejji.
\newblock Solving a system of nonlinear functional equations using revised new
  iterative method.
\newblock {\em International Journal of Mathematical and Computational
  Sciences}, 6(8):968--972, 2012.

\bibitem{bocharov2000numerical}
G.~A. Bocharov and F.~A. Rihan.
\newblock Numerical modelling in biosciences using delay differential
  equations.
\newblock {\em Journal of Computational and Applied Mathematics},
  125(1-2):183--199, 2000.

\bibitem{bui2010explicit}
T.~Bui.
\newblock Explicit and implicit methods in solving differential equations.
\newblock 2010.

\bibitem{butcher2000numerical}
J.~C. Butcher.
\newblock Numerical methods for ordinary differential equations in the 20th
  century.
\newblock {\em Journal of Computational and Applied Mathematics},
  125(1-2):1--29, 2000.

\bibitem{daftardar2008solving}
V.~Daftardar-Gejji and S.~Bhalekar.
\newblock Solving fractional diffusion-wave equations using a new iterative
  method.
\newblock {\em Fractional Calculus and Applied Analysis}, 11(2):193p--202p,
  2008.

\bibitem{daftardar2010solving}
V.~Daftardar-Gejji and S.~Bhalekar.
\newblock Solving fractional boundary value problems with dirichlet boundary
  conditions using a new iterative method.
\newblock {\em Computers \& Mathematics with Applications}, 59(5):1801--1809,
  2010.

\bibitem{daftardar2006iterative}
V.~Daftardar-Gejji and H.~Jafari.
\newblock An iterative method for solving nonlinear functional equations.
\newblock {\em Journal of Mathematical Analysis and Applications},
  316(2):753--763, 2006.

\bibitem{daftardar2015solving}
V.~Daftardar-Gejji, Y.~Sukale, and S.~Bhalekar.
\newblock Solving fractional delay differential equations: A new approach.
\newblock {\em Fractional Calculus and Applied Analysis}, 18(2):400--418, 2015.

\bibitem{dehestani2020numerical}
H.~Dehestani, Y.~Ordokhani, and M.~Razzaghi.
\newblock Numerical technique for solving fractional generalized
  pantograph-delay differential equations by using fractional-order hybrid
  bessel functions.
\newblock {\em International Journal of Applied and Computational Mathematics},
  6(1):1--27, 2020.

\bibitem{enright1995interpolating}
W.~H. Enright and M.~Hu.
\newblock Interpolating runge-kutta methods for vanishing delay differential
  equations.
\newblock {\em Computing}, 55(3):223--236, 1995.

\bibitem{evans2005adomian}
D.~J. Evans and K.~Raslan.
\newblock The adomian decomposition method for solving delay differential
  equation.
\newblock {\em International Journal of Computer Mathematics}, 82(1):49--54,
  2005.

\bibitem{ibrahim2018chaotic}
K.~Ibrahim, R.~Jamal, and F.~Ali.
\newblock Chaotic behaviour of the rossler model and its analysis by using
  bifurcations of limit cycles and chaotic attractors.
\newblock In {\em J. Phys. Conf. Ser}, volume 1003, page 012099, 2018.

\bibitem{ishak2008two}
F.~Ishak, M.~Suleiman, and Z.~Omar.
\newblock Two-point predictor-corrector block method for solving delay
  differential equations.
\newblock {\em MATEMATIKA: Malaysian Journal of Industrial and Applied
  Mathematics}, 24:131--140, 2008.

\bibitem{jhinga2019new}
A.~Jhinga and V.~Daftardar-Gejji.
\newblock A new numerical method for solving fractional delay differential
  equations.
\newblock {\em Computational and Applied Mathematics}, 38(4):1--18, 2019.

\bibitem{karoui1995numerical}
A.~Karoui and R.~Vaillancourt.
\newblock A numerical method for vanishing-lag delay differential equations.
\newblock {\em Applied Numerical Mathematics}, 17(4):383--395, 1995.

\bibitem{rihan2014delay}
F.~A. Rihan, D.~Abdelrahman, F.~Al-Maskari, F.~Ibrahim, and M.~A. Abdeen.
\newblock Delay differential model for tumour-immune response with
  chemoimmunotherapy and optimal control.
\newblock {\em Computational and mathematical methods in medicine}, 2014, 2014.

\bibitem{shampine2000solving}
L.~F. Shampine, S.~Thompson, and J.~Kierzenka.
\newblock Solving delay differential equations with dde23.
\newblock {\em URL http://www. runet. edu/\~{} thompson/webddes/tutorial. pdf},
  2000.

\bibitem{sukale2017new}
Y.~Sukale and V.~Daftardar-Gejji.
\newblock New numerical methods for solving differential equations.
\newblock {\em International Journal of Applied and Computational Mathematics},
  3(3):1639--1660, 2017.

\bibitem{suli2003introduction}
E.~S{\"u}li and D.~F. Mayers.
\newblock {\em An introduction to numerical analysis}.
\newblock Cambridge university press, 2003.

\end{thebibliography}
\bibliographystyle{abbrv}

\end{document}